\documentclass[11pt]{article}
\usepackage{amsfonts}
\usepackage{epsfig}

\def\be{\begin{equation}}
\def\ee{\end{equation}}
\def\bea{\begin{eqnarray}}
\def\eea{\end{eqnarray}}
\def\bes{\begin{eqnarray*}}
\def\ees{\end{eqnarray*}}

\def\nn{\nonumber}
\def\<{\langle}
\def\>{\rangle}
\def\lb{\label}

\def\R{{\bf R}}

\def\Z{{\bf Z}}
\def\N{{\bf N}}

\def\Q{{\bf Q}}

\def\ga{{\gamma}}

\def\th{{\theta}}

\def\Lm{{\Lambda}}

\def\rank{{\rm rank}}

\def\ol#1{\overline{#1}}  

\def\hb{\vrule height0.18cm width0.14cm $\,$}
\def\ol#1{\overline{#1}}  
\def\mapright#1{\smash{\mathop{\longrightarrow}\limits^{#1}}}

\title{The Fadell-Rabinowitz index and multiplicity of non-contractible closed geodesics on Finsler $\mathbb{R}P^{n}$}
\author{Hui Liu$^{1}$\thanks{Partially supported by NSFC (No.11401555), Anhui Provincial Natural Science Foundation (No. 1608085QA01).
E-mail: huiliu@ustc.edu.cn.} \\\\
$^{1}$ Key Laboratory of Wu Wen-Tsun Mathematics, Chinese Academy of Sciences,\\
School of Mathematical Sciences, University of Science and Technology of China,
\\Hefei, Anhui 230026, People's Republic of China\\}
\date{}

\begin{document}

\maketitle
\begin{abstract}
{\it In this paper, we prove that for every irreversible Finsler $n$-dimensional real projective
space $(\mathbb{R}P^n,F)$ with reversibility $\lambda$ and flag curvature $K$ satisfying $\frac{16}{9}\left(\frac{\lambda}{1+\lambda}\right)^2<K\le 1$
with $\lambda<3$, there exist at least $n-1$ non-contractible closed geodesics. In addition, if
the metric $F$ is bumpy with $\frac{64}{25}\left(\frac{\lambda}{1+\lambda}\right)^2<K\le 1$ and $\lambda<\frac{5}{3}$, then there
exist at least $2[\frac{n+1}{2}]$ non-contractible closed geodesics, which is the optimal lower bound due to Katok's example. The main ingredients of the proofs are
the Fadell-Rabinowitz index theory of non-contractible closed geodesics on $(\mathbb{R}P^n,F)$ and the $S^1$-equivariant Poincar$\acute{e}$ series
of the non-contractible component of the free loop space on $\mathbb{R}P^n$. }
\end{abstract}


\renewcommand{\theequation}{\thesection.\arabic{equation}}
\renewcommand{\thefigure}{\thesection.\arabic{figure}}

\setcounter{equation}{0}
\section{Introduction and main results}
In this paper, we are interested in the existence of non-contractible closed geodesics on real projective space $\mathbb{R}P^{n}$
with an irreversible Finsler metric, which is the typically non-simply connected manifold with the fundamental group $\pi_1(\mathbb{R}P^{n})=\Z_2$.

A closed curve on a Finsler manifold is a closed geodesic if it is
locally the shortest path connecting any two nearby points on this
curve. As usual, on any Finsler manifold
$(M, F)$, a closed geodesic $c:S^1=\R/\Z\to M$ is {\it prime}
if it is not a multiple covering (i.e., iteration) of any other
closed geodesics. Here the $m$-th iteration $c^m$ of $c$ is defined
by $c^m(t)=c(mt)$. The inverse curve $c^{-1}$ of $c$ is defined by
$c^{-1}(t)=c(1-t)$ for $t\in \R$.  Note that unlike Riemannian manifold,
the inverse curve $c^{-1}$ of a closed geodesic $c$
on a irreversible Finsler manifold need not be a geodesic.
We call two prime closed geodesics
$c$ and $d$ {\it distinct} if there is no $\th\in (0,1)$ such that
$c(t)=d(t+\th)$ for all $t\in\R$.
We shall omit the word {\it distinct} when we talk about more than one prime closed geodesic.
For a closed geodesic $c$ on $(M,\,F)$, denote by $P_c$
the linearized Poincar\'{e} map of $c$. Recall that a Finsler metric $F$ is {\it bumpy} if all the closed geodesics
on $(M, \,F)$ are non-degenerate, i.e., $1\notin \sigma(P_c)$ for any closed
geodesic $c$.

There have been a great deal of works on the multiplicity of closed geodesics for simply
connected manifolds (cf. \cite{Ano}, \cite{Ban}, \cite{Bot}, \cite{Fra}, \cite{Kat}, \cite{Kli}, \cite{Liu}, \cite{LiL}, \cite{Lon}, \cite{LyF}, \cite{Zil}).
In particular, Gromoll and Meyer \cite{GM2} of 1969
established the existence of infinitely many distinct closed geodesics on $M$, provided
that the Betti number sequence $\{b_p(\Lambda M; \Q)\}_{p\in\N}$ of the free loop space $\Lambda M$ of $M$ is unbounded. Then in \cite{ViS} of 1976,
for compact simply connected manifold $M$, Vigu$\acute{e}$-Poirrier and Sullivan further proved this Betti
number sequence is bounded if and only if M satisfies \bea H^*(M; \Q)\cong T_{d, n+1}(x)=\Q[x]/(x^{n+1}=0)\nn\eea
with a generator $x$ of degree $d\geq2$ and height $n+1\geq 2$, where $\dim M = dn$.
Due to their works, many mathematicians are concerned with the compact globally
symmetric spaces of rank 1 which consist in \bea S^n,\quad \mathbb{R}P^n,\quad \mathbb{C}P^n,\quad \mathbb{H}P^n \quad and \quad CaP^2,\nn\eea
where the Gromoll-Meyer assumption does not hold.

Recently the Maslov-type index theory for symplectic paths has been applied to study the closed
geodesic problem on simply connected manifolds. In 2005, Bangert and Long proved the existence of at least two distinct closed
geodesics on every Finsler $(S^2, F)$ (which was published as \cite{BaL} in 2010). Since then in the last
ten years, a great number of results on the multiplicity of closed geodesics on simply connected Finsler
manifolds have appeared, for which we refer readers to \cite{DuL1}-\cite{DuL2}, \cite{LoD}, \cite{Rad4}-\cite{Rad6},
\cite{Wan1}-\cite{Wan2}, \cite{Tai1}, \cite{HiR}, \cite{DLW1}-\cite{DLW2} and the references therein.

As for the multiplicity of closed geodesics on non-simply connected manifolds whose
free loop space possesses bounded Betti number sequence, Ballman et al. \cite{BTZ1} proved in 1981
that every Riemannian manifold with the fundamental group being a nontrivial finitely cyclic
group and possessing a generic metric has infinitely many distinct closed geodesics. In 1984,
Bangert and Hingston \cite{BaH} proved that any Riemannian manifold with the fundamental group being an
infinite cyclic group has infinitely many distinct closed geodesics. Since then, there seem to be
very few works on the multiplicity of closed geodesics on non-simply connected manifolds.
The main reason is that the topological structures of the free loop spaces on these manifolds
are not well known, so that the classical Morse theory is difficult to be applicable. In \cite{XL} of 2015,
Xiao and Long studied the topological structure of the non-contractible component of the free loop space on the real projective spaces with odd
dimensions. Furthermore, they obtain the resonance identity for the non-contractible closed
geodesics provided that the total number of closed geodesics on $\mathbb{R}P^{2n+1}$ is finite. Based on these results,
Duan, Long and Xiao \cite{DLX} further proved the existence of at least two distinct non-contractible closed
geodesics on $\mathbb{R}P^3$ endowed with a bumpy and irreversible Finsler metric. In a very recent paper \cite{Tai2}, Taimanov
studied the rational equivariant cohomology of the spaces of non-contractible loops in compact space forms and
proved the existence of at least two distinct non-contractible closed geodesics on $\mathbb{R}P^2$ endowed with a bumpy and irreversible Finsler metric.

In this paper, we continue to study the multiplicity of non-contractible closed
geodesics on $\mathbb{R}P^n$ under some natural conditions and obtain the following theorems:

{\bf Theorem 1.1.} {\it On every irreversible Finsler $n$-dimensional real projective
space  $(\mathbb{R}P^n,F)$ with reversibility $\lambda$ and flag curvature $K$ satisfying $\frac{16}{9}\left(\frac{\lambda}{1+\lambda}\right)^2<K\le 1$
with $\lambda<3$, there exist at least $n-1$ non-contractible closed geodesics.}

{\bf Theorem 1.2.} {\it On every bumpy irreversible Finsler $n$-dimensional real projective
space  $(\mathbb{R}P^n,F)$ with reversibility $\lambda$ and flag curvature $K$ satisfying $\frac{64}{25}\left(\frac{\lambda}{1+\lambda}\right)^2<K\le 1$
with $\lambda<\frac{5}{3}$, there exist at least $2[\frac{n+1}{2}]$ non-contractible closed geodesics.}

{\bf Theorem 1.3.} {\it On every irreversible Finsler $n$-dimensional real projective
space  $(\mathbb{R}P^{n},F)$ with reversibility $\lambda$ and flag curvature $K$ satisfying $\frac{4m^2}{(n-1)^2}\left(\frac{\lambda}{1+\lambda}\right)^2<K\le 1$
with $\frac{n-1}{2}\leq m<n-1$ and $\lambda<\frac{n-1}{2m-n+1}$, there exist at least $m+1$ non-contractible closed geodesics.}

{\bf Remark 1.4.} (i) For Theorem 1.1, similar results were obtained for closed geodesics on pinched Riemannian spheres and Finsler
compact simply-connected manifolds in \cite{BTZ2} and \cite{Rad4} respectively. In particular,
Ballman et al. \cite{BTZ2} applied Lusternik-Schnirelmann theory to prove that for a Riemannian metric on $S^n$ with
sectional curvature $1/4\le K\le 1$ there exist $g(n)$ geometrically distinct closed geodesics without self-intersections, and $\frac{n(n+1)}{2}$ geometrically
distinct closed geodesics without self-intersections if the metric is bumpy. Rademacher \cite{Rad4} used Fadell-Rabinowitz index in a relative version to
prove the existence of at least $n/2-1$ prime closed geodesics with length $<2n\pi$ on every Finsler  $n$-sphere
$(S^n,\,F)$ satisfying $\left(\frac{\lambda}{\lambda+1}\right)^2<K\le 1$. Differently from \cite{BTZ2} and \cite{Rad4}, we
are concerned with the number of closed geodesics on non-simply connected manifolds $\mathbb{R}P^n$
and all the closed geodesics obtained are non-contractible, our proof uses an absolute version of the
Fadell-Rabinowitz index which we restrict it to the non-contractible component of the free loop space on $\mathbb{R}P^n$.
Note that an analogous result for closed characteristics on pinched compact convex hypersurfaces in $\R^{2n}$ was given by Ekeland and Lasry \cite{EkL}.

(ii) A. Katok \cite{Kat} in 1973 found some
non-symmetric Finsler metrics on $S^n$ with only finitely many prime
closed geodesics. The smallest number of closed geodesics
on $S^n$ that one obtains in these examples is $2[\frac{n+1}{2}]$ (cf. \cite{Zil}).
Then D. V. Anosov in I.C.M. of 1974 conjectured that
the lower bound of the number of closed geodesics on
any Finsler sphere $(S^n,\,F)$ should be $2[\frac{n+1}{2}]$,
i.e., the number of closed geodesics in Katok's example. Let $(S^n, N_\alpha)$ denote the sphere $S^n$ endowed with the Katok metric,
then it induces a Finsler metric for $\mathbb{R}P^n$, which is still denoted by $N_\alpha$ for simplicity. Then
there exist precisely $2[\frac{n+1}{2}]$ non-contractible closed geodesics on $(\mathbb{R}P^n, N_\alpha)$ (cf. Proposition 3.1,
Remark 3.2 of \cite{XL}). Thus it is natural to conjecture that the lower bound of the number of non-contractible closed geodesics on
any Finsler real projective space $(\mathbb{R}P^n,\,F)$ should be $2[\frac{n+1}{2}]$.
Our Theorem 1.2 gives a confirmed answer to this conjecture for a generic case. The crucial point of the proof is the
$S^1$-equivariant Poincar$\acute{e}$ series of the non-contractible component of the free loop space on $\mathbb{R}P^n$
derived by Taimanov \cite{Tai2}.

(iii) For Theorem 1.3, a similar result for the number of closed geodesics on positively curved compact simply-connected manifolds was given in Corollary 1.3
of \cite{Rad2}. In our Theorem 1.3, for non-simply connected manifolds $\mathbb{R}P^n$, we get a lower bound for the number of closed geodesics
which are non-contractible, so our result is different from that of \cite{Rad2}.

Note that a closed geodesic on $\mathbb{R}P^n$ is orientable if and only if $n$ is odd, the iteration formulaes of Maslov-type index
of orientable and non-orientable closed geodesics are different, cf. \cite{LiL}, \cite{Liu}. In the proofs of our theorems,
we will not use the iteration formulae of Maslov-type index, but we use Morse index and rational coefficient homology only,
so the orientable or non-orientable problem is inessential for our proofs.

In this paper, let $\N$, $\N_0$, $\Z$, $\Z_2$, $\Q$, $\R$ denote the sets of natural integers,
non-negative integers, integers, cyclic group of order 2, rational numbers, real numbers respectively.
We define the function $[a]=\max{\{k\in {\bf Z}\mid k\leq a\}}$.

\setcounter{equation}{0}
\section{Morse theory of non-contractible closed geodesics}
Let $M=(M,F)$ be a compact Finsler manifold, the space
$\Lambda=\Lambda M$ of $H^1$-maps $\gamma:S^1\rightarrow M$ has a
natural structure of Riemannian Hilbert manifolds on which the
group $S^1=\R/\Z$ acts continuously by isometries (cf. \cite{Kli}, \cite{She}). This action is defined by
$(s\cdot\gamma)(t)=\gamma(t+s)$ for all $\gamma\in\Lm$ and $s,
t\in S^1$. For any $\gamma\in\Lambda$, the energy functional is
defined by
\be E(\gamma)=\frac{1}{2}\int_{S^1}F(\gamma(t),\dot{\gamma}(t))^2dt.
\lb{2.1}\ee
It is $C^{1,1}$ and invariant under the $S^1$-action. The
critical points of $E$ of positive energies are precisely the closed geodesics
$\gamma:S^1\to M$. The index form of the functional $E$ is well
defined along any closed geodesic $c$ on $M$, which we denote by
$E''(c)$. As usual, we denote by $i(c)$ and
$\nu(c)$ the Morse index and nullity of $E$ at $c$. In the
following, we denote by
\be \Lm^\kappa=\{d\in \Lm\;|\;E(d)\le\kappa\},\quad \Lm^{\kappa-}=\{d\in \Lm\;|\; E(d)<\kappa\},
  \quad \forall \kappa\ge 0. \nn\ee
For a closed geodesic $c$ we set $ \Lm(c)=\{\ga\in\Lm\mid E(\ga)<E(c)\}$.

For $m\in\N$ we denote the $m$-fold iteration map
$\phi_m:\Lambda\rightarrow\Lambda$ by $\phi_m(\ga)(t)=\ga(mt)$, for all
$\,\ga\in\Lm, t\in S^1$, as well as $\ga^m=\phi_m(\gamma)$. If $\gamma\in\Lambda$
is not constant then the multiplicity $m(\gamma)$ of $\gamma$ is the order of the
isotropy group $\{s\in S^1\mid s\cdot\gamma=\gamma\}$. For a closed geodesic $c$,
the mean index $\hat{i}(c)$ is defined as usual by
$\hat{i}(c)=\lim_{m\to\infty}i(c^m)/m$. Using singular homology with rational
coefficients we consider the following critical $\Q$-module of a closed geodesic
$c\in\Lambda$:
\be \overline{C}_*(E,c)
   = H_*\left((\Lm(c)\cup S^1\cdot c)/S^1,\Lm(c)/S^1; \Q\right). \lb{2.3}\ee
Following \cite{Rad1}, Section 6.2, we can use finite-dimensional
approximations to $\Lambda$ to apply the results of D. Gromoll and W. Meyer
\cite{GM1} to a given closed geodesic $c$ which is isolated as a critical orbit.
Then we have

{\bf Proposition 2.1.} {\it Let $k_j(c)\equiv\dim\overline{C}_j(E,c)$.
Then $k_j(c)$  equal to $0$ when $j<i(c)$ or $j>i(c)+\nu(c)$ and if $\nu(c)=0$, then $k_j(c)=1$ for $j=i(c)$. }

For $M=\mathbb{R}P^{n}$, it is well known that $\pi_1(\mathbb{R}P^{n})=\Z_2=\{e, g\}$ with $e$ being
the identity and $g$ being the generator of $\Z_2$ satisfying $g^2=e$. Then the free loop space $\Lambda M$ possesses
a natural decomposition\bea \Lambda M=\Lambda_e M\bigsqcup\Lambda_g M,\nn\eea
where $\Lambda_e M$ and $\Lambda_g M$ are the two connected components of $\Lambda M$ whose elements are homotopic
to $e$ and $g$ respectively.

{\bf Lemma 2.2.} (cf. Theorem 3 of \cite{Tai2}) {\it For $M=\mathbb{R}P^n$, we have

(i) When $n=2k+1$ is odd, the $S^1$-cohomology ring of $\Lambda_g M$ has the form
$$H^{S^1, *}(\Lambda_g M; \Q)=\Q[w, z]/ \{w^{k+1} = 0\}, \quad deg(w)=2, \quad deg(z)=2k$$
Then the $S^1$-equivariant Poincar$\acute{e}$ series
of $\Lambda_g M$ is given by
\bea P^{S^1}(\Lambda_g M; \Q)(t)&=&\frac{1-t^{2k+2}}{(1-t^2)(1-t^{2k})}\nn\\
&=&\frac{1}{1-t^2}+\frac{t^{2k}}{1-t^{2k}}\nn\\&=&(1+t^2+t^4+\cdots+t^{2k}+\cdots)+(t^{2k}+t^{4k}+t^{6k}+\cdots),\nn\eea
which yields Betti numbers
\bea b_q
&=& \rank H_q^{S^1}(\Lambda_g M;\Q)  \nn\\
&=& \left\{\matrix{
    2,&\quad {\it if}\quad q\in \{j(n-1)\mid j\in\N\},  \cr
    1,&\quad {\it if}\quad q\in2\N_0\setminus \{j(n-1)\mid j\in\N\},  \cr
    0, &\quad {\it otherwise}. \cr}\right. \lb{2.4}\eea

(ii) When $n=2k$ is even, the $S^1$-cohomology ring of $\Lambda_g M$ has the form
$$H^{S^1, *}(\Lambda_g M; \Q)=\Q[w, z]/ \{w^{2k} = 0\}, \quad deg(w)=2, \quad deg(z)=4k-2$$
Then the $S^1$-equivariant Poincar$\acute{e}$ series
of $\Lambda_g M$ is given by \bea P^{S^1}(\Lambda_g M; \Q)(t)&=&\frac{1-t^{4k}}{(1-t^2)(1-t^{4k-2})}\nn\\
&=&\frac{1}{1-t^2}+\frac{t^{4k-2}}{1-t^{4k-2}}\nn\\&=&(1+t^2+t^4+\cdots+t^{2k}+\cdots)+(t^{4k-2}+t^{2(4k-2)}+t^{3(4k-2)}+\cdots),\nn\eea
which yields Betti numbers \bea b_q
&=& \rank H_q^{S^1}(\Lambda_g M;\Q)  \nn\\
&=& \left\{\matrix{
    2,&\quad {\it if}\quad q\in \{2j(n-1)\mid j\in\N\},  \cr
    1,&\quad {\it if}\quad q\in2\N_0\setminus \{2j(n-1)\mid j\in\N\},  \cr
    0, &\quad {\it otherwise}. \cr}\right. \lb{2.5}\eea}

Note that for a non-contractible prime closed geodesic $c$, $c^m\in\Lambda_g M$ if and only if $m$ is odd. Thus
we have the following Morse inequality for non-contractible closed geodesics:

{\bf Theorem 2.3.} (cf. Theorem I.4.3 of \cite{Cha})
{\it Assume that $M=\mathbb{R}P^n$ be a Finsler manifold with finitely many non-contractible prime closed geodesics, denoted by $\{c_j\}_{1\le j\le k}$. Set
\bea M_q =\sum_{1\le j\le k,\; m\ge 1}\dim{\ol{C}}_q(E, c^{2m-1}_j),\quad q\in\Z.\nn\eea
Then for every integer $q\ge 0$ there holds }
\bea M_q - M_{q-1} + \cdots +(-1)^{q}M_0
&\ge& {b}_q - {b}_{q-1}+ \cdots + (-1)^{q}{b}_0, \nn\\
M_q &\ge& {b}_q. \lb{2.6}\eea

\setcounter{equation}{0}
\section{Fadell-Rabinowitz index theory for non-contractible closed geodesics on $(\mathbb{R}P^{n}, F)$}
In 1978, Fadell and Rabinowitz \cite{FaR} introduced a topological index theory by using the equivariant Borel
cohomology. Via a classifying map the cohomology ring $H^*(BG)$ of a classifying space is a subring of the equivariant cohomology
$H_G^*(X)$ of the G-space $X$. For a characteristic class $\eta\in H^*(BG)$ the Fadell-Rabinowitz index $index_\eta X$ of X
is the smallest non-negative integer k with $\eta^k=0$ in $H_G^*(X)$. Hence the classes $\eta^j$
for $j=0, \ldots, index_\eta X- 1$ define a sequence of subordinate cohomology classes.
Fadell and Rabinowitz used this index to study the bifurcation of time periodic solutions
from an equilibrium solution for Hamiltonian systems of ordinary differential
equations. Ekeland and Hofer in \cite{EkH} of 1987 used this index for the group $S^1$ to
derive a close relationship between the set of Maslov-type indices of closed characteristics on convex
Hamiltonian energy surfaces in $\R^{2n}$ and the set of even
positive integers, which is the core in studying multiplicity and ellipticity of
closed characteristics on compact convex hypersurfaces (cf. \cite{LoZ}). Rademacher in \cite{Rad2} of 1994 developped
the relative version of Fadell-Rabinowitz index to study the existence of closed geodesics on compact simply connected
manifolds, which played an important role in studying multiplicity and stability of
closed geodesics on Finsler spheres, cf. \cite{Rad4}, \cite{Wan3}.

In this section, we use only singular homology modules with $\Q$-coefficients and apply the Fadell-Rabinowitz
index to the non-contractible component $\Lambda_g M$ of the free
loop space on $M=\mathbb{R}P^{n}$, which is one of the main ingredients for proving our multiplicity results of non-contractible closed geodesics.
The main idea is to extend the work of Rademacher \cite{Rad2} to the non-contractible component $\Lambda_g M$,
thus the following lemmas and theorems look like those of \cite{Rad2}.

It is well known that there exists a non-contractible closed geodesic $c$ on $(\mathbb{R}P^n, F)$, which
is the minimal point of the energy functional $E$ on $\Lambda_g$, set $E(c)=\frac{\kappa_1^2}{2}$ for $\kappa_1>0$.
For any $\kappa>0$, we denote by
\be \Lambda_g^\kappa=\{\gamma\in \Lambda_g\;|\;
             E(\gamma)\leq\frac{1}{2} \kappa^2\}. \lb{3.1}\ee
We define the function
$d_g: [\kappa_1, +\infty]\rightarrow\N\cup\{\infty\}$ by \bea d_g(\kappa)=index_\eta \Lambda_g^\kappa,\nn\eea
where $H^*(BS^1)=H^*(CP^\infty)=\Q[\eta]$.

{\bf Definition 3.1.} {\it The lower index $\underline{\sigma}$ and upper index $\overline{\sigma}$ are given by:
\bea\underline{\sigma}=\lim\inf_{\kappa\rightarrow\infty}{\frac{d_g(\kappa)}{\kappa}},\quad
\overline{\sigma}=\lim\sup_{\kappa\rightarrow\infty}{\frac{d_g(\kappa)}{\kappa}},\nn\eea
we call $[\underline{\sigma}, \overline{\sigma}]$ the global index interval.}

Similarly to Proposition 5.2 of \cite{Rad2}, as a consequence of the Morse-Schoenberg comparison theorem and the
proof of Myers' theorem, we have the following estimates for the global index interval
in terms of the curvature.

{\bf Lemma 3.2.} {\it Let $K$ and $Ric$ be the flag curvature and Ricci curvature  of the manifold $(\mathbb{R}P^n, F)$ respectively, then we have

(a) If $K\leq \Delta^2$, then  $\overline{\sigma}\leq \frac{\Delta(n-1)}{2\pi}$.

(b) If $K\geq \delta^2$, then  $\overline{\sigma}\geq \frac{\delta(n-1)}{2\pi}$.

(c) If $Ric\geq \delta^2(n-1)$, then  $\underline{\sigma}\geq \frac{\delta}{2\pi}$.}

{\bf Proof.} The proof is somewhat similar to Proposition 5.2 of \cite{Rad2},
for the reader's convenience, we sketch a proof here. We only prove (a).
Choose a sequence $(F_i)$ of bumpy metrics converging to $F$ in the strong
$C^2$ topology. Let $E_i(\gamma)=\frac{1}{2}\int_{S^1}F_i(\dot{\gamma}(t))^2dt$
be the energy functional and $\Lambda_i^\kappa=\{\gamma\in \Lambda_g\;|\; E_i(\gamma)\leq\frac{1}{2} \kappa^2\}$.

Choose $(a_i)\subset\R$ with $a=\lim_{i\rightarrow\infty} a_i$ and $\Lambda_g^a\subset\Lambda_i^{a_i}$ for all $i$. We can
choose a sequence $(\Delta_i)\subset\R$ such that the flag curvature of $F_i$ is bounded from
above by $\Delta_i^2$ and $\Delta=\lim_{i\rightarrow\infty} \Delta_i$. It follows from the Morse-Schoenberg comparison
theorem that the index $i(c)$ of a closed geodesic $c$ of the metric $F_i$ with length at most $a_i$ satisfies
\bea i(c)\leq \left(\frac{a_i}{\pi}\Delta_i+1\right)(n-1)=:k_i.\lb{3.2}\eea
Hence it follows from the Morse lemma that $H^k_{S^1}(\Lambda_i^{a_i})=0$ for all $k>k_i$.
Since $\Lambda_g^a\subset\Lambda_i^{a_i}$ it follows from the composition
\bea H^k(BS^1)\rightarrow H_{S^1}^k(\Lambda_i^{a_i})\rightarrow  H_{S^1}^k(\Lambda_g^{a}).\nn\eea
of restriction homomorphisms that $\eta^k=0$ in $H_{S^1}^{2k}(\Lambda_g^{a})$ for all $k>\frac{k_i}{2}$.
Hence we have $d_g(a)\leq \frac{k_i}{2}$, which together with (\ref{3.2}) and Definition 3.1 gives
$\overline{\sigma}\leq \frac{\Delta(n-1)}{2\pi}$.\hfill\hb

{\bf Corollary 3.3.} {\it For $\mathbb{R}P^n$ endowed with any Finsler metric $F$, we have $\lim_{\kappa\rightarrow+\infty}d_g(\kappa)=+\infty$.}

{\bf Proof.} Let $(S^n, N_\alpha)$ denote the sphere $S^n$ endowed with the Katok metric, then it induces
a Finsler metric for $\mathbb{R}P^n$, which is still denoted by $N_\alpha$ for simplicity. Therefore the natural projection
\bea \pi:(S^n, N_\alpha)\rightarrow (\mathbb{R}P^n, N_\alpha)\nn\eea
is locally isometric. Since $(S^n, N_\alpha)$ and $(\mathbb{R}P^n, N_\alpha)$ have constant flag curvature 1, then for $(\mathbb{R}P^n, N_\alpha)$,
$\underline{\sigma}>0$ by Lemma 3.2. Thus for $\mathbb{R}P^n$ endowed with any Finsler metric $F$, $\underline{\sigma}>0$
also holds which implies $\lim_{\kappa\rightarrow+\infty}d_g(\kappa)=+\infty$. In fact, Let $F,F^*$ be two Finsler metrics on $\mathbb{R}P^n$ with
\bea F^*(X)^2/D^2\leq  F(X)^2\leq D^2F^*(X)^2\nn\eea
for all tangent vector $X$ and some $D>1$. Let $[\underline{\sigma}, \overline{\sigma}]$, $[\underline{\sigma}^*, \overline{\sigma}^*]$
be the global index intervals with respect to the metrics $F,F^*$ respectively. Then by definition we have
\bea D^{-1}\underline{\sigma}\leq \underline{\sigma}^*\leq D\underline{\sigma},
\quad D^{-1}\overline{\sigma}\leq\overline{\sigma}^*\leq D\overline{\sigma}\nn\eea
which implies our desired result.\hfill\hb

{\bf Lemma 3.4.} {\it The function $d_g:  [\kappa_1, +\infty)\rightarrow\N$ is non-decreasing and
$\lim_{\lambda\searrow \kappa}d_g(\lambda)=d_g(\kappa)$.
Each discontinuous point of $d_g$ is a critical value of the energy functional $E$ on $\Lambda_g$.
In particular, if $d_g(\kappa)-d_g(\kappa-)\ge 2$, then there are infinitely
many prime non-contractible closed geodesics $c$ with energy $\kappa$, where $d_g(\kappa-)=\lim_{\epsilon\searrow 0}d_g(\kappa-\epsilon)$,
$t\searrow a$  means $t>a$ and $t\to a$.}

{\bf Proof.} We have that $d_g(\lambda)$ is finite by the proof of Lemma 3.2.
The rest of our proof is completely similar to Lemmas 5.5-5.6 and Corollary 5.7 of \cite{Rad2}, thus
we omit it.\hfill\hb

For each $i\ge 1$, define
\bea \kappa_i=\inf\{\delta\in\R\;|\: d_g(\delta)\ge i\}.\lb{3.3}\eea
Then $\kappa_i$ is well defined by Corollary 3.3.

Now we can use the methods of Lemmas 2.3-2.4 of \cite{Wan3} to obtain the main results in this section.

{\bf Theorem 3.5.} {\it Suppose there are only finitely many prime non-contractible closed
geodesics on $(\mathbb{R}P^n,\, F)$. Then each $\kappa_i$ is a critical value of $E$ on $\Lambda_g$.
If $\kappa_i=\kappa_j$ for some $i<j$, then there are infinitely many prime non-contractible
closed geodesics on $(\mathbb{R}P^n,\, F)$. }

{\bf Proof.} We replace Lemma 2.2 used in the proof of Lemma 2.3 of \cite{Wan3} by the above Lemma 3.4, then
our proof is similar to that of Lemma 2.3 of \cite{Wan3}.\hfill\hb

{\bf Theorem 3.6.} {\it Suppose there are only finitely many prime non-contractible closed
geodesics on $(\mathbb{R}P^n,\, F)$. Then for every $i\in\N$, there exists a non-contractible
closed geodesic $c$ on $(\mathbb{R}P^n,\, F)$ such that} \bea E(c)=\kappa_i,\quad
\ol{C}_{2i-2}(E, c)\neq 0.\lb{3.4}\eea

{\bf Proof.} Choose $\epsilon$ small enough such that the interval
$(\kappa_i-\epsilon,\,\kappa_i+\epsilon)$ contains no critical values of $E$ except $\kappa_i$.
By the same proof of (2.12) of \cite{Wan3}, we have
\be H^{\ast}_{S^1}(\Lambda_g^{\kappa_i+\epsilon},\;\Lambda_g^{\kappa_i-\epsilon})
\cong H^\ast(\Lambda_g^{\kappa_i+\epsilon}/S^1,\;\Lambda_g^{\kappa_i-\epsilon}/S^1).\lb{3.5}\ee
Also as (2.13) of \cite{Wan3}, we have
\be H_\ast(\Lambda_g^{\kappa_i+\epsilon}/S^1,\;\Lambda_g^{\kappa_i-\epsilon}/S^1)
\cong H^\ast(\Lambda_g^{\kappa_i+\epsilon}/S^1,\;\Lambda_g^{\kappa_i-\epsilon}/S^1).\lb{3.6}\ee
Then we combining (\ref{3.5})-(\ref{3.6}) with Theorem 1.4.2 of \cite{Cha} to obtain
\be H_{S^1}^*(\Lambda_g^{\kappa_i+\epsilon},\;\Lambda_g^{\kappa_i-\epsilon})
=\bigoplus_{E(c)=\kappa_i}\ol{C}_{\ast}(E, \; c).\lb{3.7}\ee
For an $S^1$-space $X$, we denote by
$X_{S^1}$ the homotopy quotient of $X$ by $S^1$, i.e., $X_{S^1}=S^\infty\times_{S^1}X$,
where $S^\infty$ is the unit sphere in an infinite dimensional complex Hilbert space.
Similar to P.431 of \cite{EkH}, we have
\be H^{2(i-1)}((\Lambda_g^{\kappa_i+\epsilon})_{S^1},\,(\Lambda_g^{\kappa_i-\epsilon})_{S^1})
\mapright{q^\ast} H^{2(i-1)}((\Lambda_g^{\kappa_i+\epsilon})_{S^1} )
\mapright{p^\ast}H^{2(i-1)}((\Lambda_g^{\kappa_i-\epsilon})_{S^1}),\lb{3.8}\ee
where $p$ and $q$ are natural inclusions. Denote by
$f: (\Lambda_g^{\kappa_i+\epsilon})_{S^1}\rightarrow CP^\infty$ a classifying map and let
$f^{\pm}=f|_{(\Lambda_g^{\kappa_i\pm\epsilon})_{S^1}}$. Then clearly each
$f^{\pm}: (\Lambda_g^{\kappa_i\pm\epsilon})_{S^1}\rightarrow CP^\infty$ is a classifying
map on $(\Lambda_g^{\kappa_i\pm\epsilon})_{S^1}$. Let $\eta \in H^2(CP^\infty)$ be the first
universal Chern class.

By definition of $\kappa_i$, we have $d_g(\kappa_i-\epsilon)< i$, hence
$(f^-)^\ast(\eta^{i-1})=0$. Note that
$p^\ast(f^+)^\ast(\eta^{i-1})=(f^-)^\ast(\eta^{i-1})$.
Hence the exactness of (\ref{3.8}) yields a
$\sigma\in H^{2(i-1)}((\Lambda_g^{\kappa_i+\epsilon})_{S^1},\,(\Lambda_g^{\kappa_i-\epsilon})_{S^1})$
such that $q^\ast(\sigma)=(f^+)^\ast(\eta^{i-1})$.
Since $d_g(\kappa_i+\epsilon)\ge i$, we have $(f^+)^\ast(\eta^{i-1})\neq 0$.
Hence $\sigma\neq 0$, and then
$$H_{S^1}^{2i-2}(\Lambda_g^{\kappa_i+\epsilon},\;\Lambda_g^{\kappa_i-\epsilon})=
H^{2i-2}((\Lambda_g^{\kappa_i+\epsilon})_{S^1},\,(\Lambda_g^{\kappa_i-\epsilon})_{S^1})\neq 0. $$
Thus (\ref{3.4}) follows from (\ref{3.7}).\hfill\hb

Note that by Proposition 2.1, (\ref{3.4}) yields
\bea E(c)=\kappa_i,\quad
i(c)\leq 2i-2 \leq i(c)+\nu(c).\lb{3.9}\eea

Similar to Definition 1.4 of \cite{LoZ}, we have

{\bf Definition 3.7.} {\it A prime non-contractible closed geodesic $c$ is
$(2m-1, i)$-{\bf variationally visible}, if there exist some $m, i\in\N$
such that (\ref{3.4}) holds for $c^{2m-1}$ and $\kappa_i$.
Denote by $\mathcal{V}(\mathbb{R}P^n, F)$ the set of
variationally visible non-contractible closed geodesics on $(\mathbb{R}P^n,F)$.}

Similar to Theorem V.3.15 of \cite{Eke} and Theorem 6.4 of \cite{Rad2}, we have

{\bf Theorem 3.8.} {\it  Assume that $\,^{\#}\mathcal{V}(\mathbb{R}P^n, F)<+\infty$
and $\hat{i}(c)>0$ for any $c\in\mathcal{V}(\mathbb{R}P^n, F)$, then we have}
\bea \sum_{c\in\mathcal{V}(\mathbb{R}P^n, F)}\frac{1}{\hat{i}(c)}\geq 1.\nn\eea

{\bf Proof.} Noticing that only odd iterations of prime non-contractible closed geodesics
can be critical points of the energy functional $E$ on $\Lambda_g \mathbb{R}P^n$,  then our proof follows from that of Theorem V.3.15 of \cite{Eke}, if we replace
Proposition 7 used in the proof of Theorem V.3.15 of \cite{Eke} by the above Theorem 3.6 and (\ref{3.9}).\hfill\hb

\setcounter{equation}{0}
\section{Multiplicity of non-contractible closed geodesics on $(\mathbb{R}P^{n}, F)$}
In this section, we give proofs of the multiplicity results of non-contractible closed geodesics on $(\mathbb{R}P^{n}, F)$
stated in Section 1 by using Fadell-Rabinowitz theory for non-contractible closed geodesics obtained in Section 3,
the $S^1$-equivariant Poincar$\acute{e}$ series of the non-contractible component of the free loop space on $\mathbb{R}P^n$
and some estimations on the index and mean index of non-contractible closed geodesics.

The following lemmas proved by Rademacher are useful for us to control the index of non-contractible closed geodesics.

{\bf Lemma 4.1.} (cf. Lemma 1 of \cite{Rad4}) {\it  Let $c$ be a closed geodesic on a Finsler manifold $(M, F)$
of dimension $n$ with positive flag curvature $K\geq\delta$ for some $\delta>0$,
if the length $L(c)$ satisfies $L(c)>\frac{k\pi}{\sqrt{\delta}}$ for some positive integer $k$ then $i(c)\geq k(n-1)$.}

{\bf Lemma 4.2.} (cf. Theorem 1 and Theorem 4 of \cite{Rad3}) {\it Let $M$ be a compact and simply-connected manifold
with a Finsler metric $F$ with reversibility $\lambda$ and flag curvature $K$ satisfying
$0< K \leq 1$ resp. $\frac{\lambda^2}{(\lambda+1)^2} < K\leq 1$ if the dimension $n$ is odd. Then the length of a
closed geodesic is bounded from below by $\pi \frac{\lambda+1}{\lambda}$.}

Now we can estimate the index of non-contractible closed geodesics.

{\bf Proposition 4.3.} {\it Let $c$ be a non-contractible closed geodesic on Finsler $n$-dimensional real projective
space  $(\mathbb{R}P^n,F)$ with reversibility $\lambda$ and flag curvature $K$ satisfying $\frac{16}{9}\left(\frac{\lambda}{1+\lambda}\right)^2<K\le 1$
with $\lambda<3$, then $i(c^m)\geq 2(n-1)$ for $m\geq3$.}

{\bf Proof.} Since $(S^n, F)$ is the universal double covering of $(\mathbb{R}P^n, F)$, then $c^2$ is a closed geodesic on $(S^n, F)$.
By Lemma 4.2, for $m\geq3$ we have
\bea L(c^m)=\frac{m}{2}L(c^{2})\geq \frac{m\pi}{2} \frac{\lambda+1}{\lambda}\geq \frac{3\pi}{2} \frac{\lambda+1}{\lambda}=
2\frac{3\pi(\lambda+1)}{4\lambda}.\lb{4.1}\eea
By assumption, we can choose an $\delta>\frac{16}{9}\left(\frac{\lambda}{1+\lambda}\right)^2$ such that $K\geq \delta$, it together
with (\ref{4.1}) gives  $L(c^m)> \frac{2\pi}{\sqrt{\delta}}$ for $m\geq3$ which implies $i(c^m)\geq 2(n-1)$ by Lemma 4.1.\hfill\hb

{\bf Proposition 4.4.} {\it Let $c$ be a non-contractible closed geodesic on Finsler $n$-dimensional real projective
space  $(\mathbb{R}P^n,F)$ with reversibility $\lambda$ and flag curvature $K$ satisfying $\frac{64}{25}\left(\frac{\lambda}{1+\lambda}\right)^2<K\le 1$
with $\lambda<\frac{5}{3}$, then $i(c^m)\geq 4(n-1)$ for $m\geq5$.}

{\bf Proof.} Since $(S^n, F)$ is the universal double covering of $(\mathbb{R}P^n, F)$, then $c^2$ is a closed geodesic on $(S^n, F)$.
By Lemma 4.2, for $m\geq5$ we have
\bea L(c^m)=\frac{m}{2}L(c^{2})\geq \frac{m\pi}{2} \frac{\lambda+1}{\lambda}\geq \frac{5\pi}{2} \frac{\lambda+1}{\lambda}=
4\frac{5\pi(\lambda+1)}{8\lambda}.\lb{4.2}\eea
By assumption, we can choose an $\delta>\frac{64}{25}\left(\frac{\lambda}{1+\lambda}\right)^2$ such that $K\geq \delta$, it together
with (\ref{4.2}) gives  $L(c^m)> \frac{4\pi}{\sqrt{\delta}}$ for $m\geq5$ which implies $i(c^m)\geq 4(n-1)$ by Lemma 4.1.\hfill\hb

{\bf Proof of Theorem 1.1.} We assume that there exist finitely many prime non-contractible closed
geodesics. Then by Theorems 3.5-3.6, there exist non-contractible
closed geodesics $c_i$ such that \bea E(c_i)=\kappa_i,\quad
\ol{C}_{2i-2}(E, c_i)\neq 0,\lb{4.3}\eea
and $\kappa_i\neq\kappa_j$ for positive integers $i\neq j$. Combining Proposition 2.1 with (\ref{4.3}),
we obtain $i(c_i)\leq2i-2\leq i(c_i)+\nu(c_i)$ for any $i\in\N$. Then by Proposition 4.3, note that
even iterations of $c_i$ do not belong to the non-contractible component $\Lambda_g$ of the free
loop space, we obtain that $c_i$ is prime for $1\leq i\leq n-1$. Then $c_i\neq c_j$ for $1\leq i\neq j\leq n-1$
since $ E(c_i)=\kappa_i\neq\kappa_j= E(c_j)$. Thus we get
$n-1$ non-contractible prime closed geodesics $\{c_i\mid 1\leq i\leq n-1\}$.\hfill\hb

{\bf Proof of Theorem 1.2.} Denote all the prime non-contractible closed
geodesics by $\{c_i\mid 1\leq i\leq k\}$. Then by Proposition 4.4, we
have $i(c_i^m)\geq 4(n-1)$ for $m\geq 5$ which implies that if $c_i^m$ has contribution to
$m_q$ for $0\leq q\leq 4(n-1)-2$, then $m=1$ or $3$ since even iterations of $c_i$ have no contribution to
$m_q$. Thus we obtain
\bea  m_{4n-6}+m_{4n-8}+\cdots+m_2+m_0\leq 2k.\lb{4.4}\eea
In the following, we prove in two cases:

(i) $n$ is odd. In this case, by (\ref{2.6}) of Theorem 2.3 and (\ref{2.4}) of Lemma 2.2, we have
\bea  m_{4n-6}+m_{4n-8}+\cdots+m_2+m_0\geq b_{4n-6}+b_{4n-8}+\cdots+b_2+b_0=2n+1.\lb{4.5}\eea
Then $k\geq n+1$ follows from (\ref{4.4}) and (\ref{4.5}).

(ii) $n$ is even. In this case, by (\ref{2.6}) of Theorem 2.3 and (\ref{2.5}) of Lemma 2.2, we have
\bea  m_{4n-6}+m_{4n-8}+\cdots+m_2+m_0\geq b_{4n-6}+b_{4n-8}+\cdots+b_2+b_0=2n-1.\lb{4.6}\eea
Then $k\geq n$ follows from (\ref{4.4}) and (\ref{4.6}). The proof is complete.\hfill\hb

The following lemma helps us to control the mean index of a closed geodesic:

{\bf Lemma 4.5.} (cf. Lemma 2 of \cite{Rad4}) {\it Let $c$ be a closed geodesic on a compact and simply-connected Riemannian
manifold of dimension $n$ with a non-reversible Finsler metric with reversibility
$\lambda$ and flag curvature $0<\delta\leq K\leq 1$ where $\delta>\frac{\lambda^2}{(\lambda+1)^2}$ if n is odd. Then
$\hat{i}(c)\geq \sqrt{\delta}\frac{\lambda+1}{\lambda}(n-1)$.}

{\bf Proof of Theorem 1.3.} Denote all the prime non-contractible closed
geodesics by $\{c_i\mid 1\leq i\leq k\}$. Since $(S^n, F)$ is the universal double covering of $(\mathbb{R}P^n, F)$,
then $c_i^2$ is a closed geodesic on $(S^n, F)$. By Lemma 4.5, we have
\bea\hat{i}(c_i^2)> \sqrt{\frac{4m^2}{(n-1)^2}\left(\frac{\lambda}{1+\lambda}\right)^2}\frac{\lambda+1}{\lambda}(n-1)=2m,\nn\eea
which implies \bea \hat{i}(c_i)>m.\lb{4.7}\eea
On the other hand, by Theorem 3.8, we have \bea \sum_{1\leq i\leq k}\frac{1}{\hat{i}(c_i)}\geq 1,\nn\eea
which together with (\ref{4.7}) gives $k\geq m+1$.  \hfill\hb

{\bf Acknowledgements.}  I would like to sincerely thank my
advisor, Professor Yiming Long, for introducing me to the theory of
closed geodesics and for his valuable advices and comments on this paper.

\bibliographystyle{abbrv}

\end{document}